\documentclass[12pt,leqno]{amsart}

\usepackage{amsmath,amsthm}
\usepackage{amssymb}
\usepackage{euscript}

\numberwithin{equation}{section}  

\newcommand{\sqsp}{\renewcommand{\baselinestretch}{1.4}\tiny\normalsize}

\newtheorem{thm}{Theorem}
\newtheorem{prop}[equation]{Proposition}

\newtheorem{cor}[equation]{Corollary}

\theoremstyle{definition}
\newtheorem{definition}[equation]{Definition}

\newcommand{\bC}{\mathbf{C}}

\newcommand{\bZ}{\mathbf{Z}}
\newcommand{\Set}{\mathbf{S}\mathbf{e}\mathbf{t}}  
\newcommand{\FRing}{\mathbf{F}\mathbf{R}\mathbf{i}\mathbf{n}\mathbf{g}} 
\newcommand{\FAlg}{\mathbf{F}\mathbf{A}\mathbf{l}\mathbf{g}}  

\DeclareMathOperator{\Id}{Id}

\newcommand{\ctensor}{\widehat{\otimes}}           
\newcommand{\llbrack}{\lbrack \lbrack}             
\newcommand{\rrbrack}{\rbrack \rbrack}             

\begin{document}
\title{Unstable $K$-cohomology algebra is filtered lambda-ring}
\author{Donald Yau}

\begin{abstract}
Boardman, Johnson, and Wilson gave a precise formulation for an unstable algebra over a generalized cohomology theory.  Modifying their definition slightly in the case of complex $K$-theory by taking into account its periodicity, we prove that an unstable algebra for complex $K$-theory is precisely a filtered $\lambda$-ring, and vice versa.
\end{abstract}

\email{dyau@math.uiuc.edu}
\address{Department of Mathematics, University of Illinois at Urbana-Champaign, 1409 W.\ Green Street, Urbana, IL 61801, USA}
\keywords{Unstable algebra, lambda-ring}
\subjclass{55N20,55N15,55S05,55S25}
\date{\today}

\maketitle

\sqsp


\section{Introduction}
\label{sec:introduction}

Lambda operations in complex $K$-theory were first introduced by Grothendieck.  These operations should be thought of as exterior power operations; in fact, for an element $\alpha$ in $K(X)$ which comes from an actual vector bundle on $X$ with $X$ a finite complex, $\lambda^i(\alpha)$ is the element represented by the $i$th exterior power of that vector bundle.  It was realized early on that these operations generate all $K$-theory operations (see \cite{toda}, for example).  A $\lambda$-ring is, roughly speaking, a commutative ring with operations $\lambda^i$ which behave exactly like $\lambda$-operations on $K$-theory of spaces.  It is, therefore, natural to think that $\lambda$-rings have all the algebraic structures to capture the unstable $K$-theory algebra of spaces.  There is, however, a more precise notion of an unstable algebra, which we now recall.

In their seminal article \cite{bjw} Boardman, Johnson, and Wilson gave a precise definition for an unstable $E^*$-cohomology algebra, where $E^*$ is a generalized cohomology theory, such as $K$-theory, satisfying some reasonable freeness conditions.  Given $E^*$ denote by $\underline{E}_k$ the $k$th space in the $\Omega$-spectrum representing $E^*$.  Of course, operations on $E^k(-)$ are just the elements of $E^*\underline{E}_k$.  There are functors   
   \[
   U^k(-) ~=~ \FAlg(E^*\underline{E}_k,-) \colon \FAlg ~\to~ \Set
   \]
from the category $\FAlg$ of complete Hausdorff filtered $E^*$-algebras and continuous $E^*$-algebra homomorphisms to sets.  Thanks to the extra structures on the spaces $\underline{E}_k$, the functor $U$ whose components are $U^k$ to graded sets becomes a comonad on the category $\FAlg$.  Then these authors define an \emph{unstable $E$-cohomology algebra} to be a $U$-coalgebra.  This definition applies to $K$-theory in particular.

The purpose of this note is to show that these two notions of unstable algebras for $K$-theory (almost) coincide.  This is, perhaps, not surprising and is even intuitively obvious.  But the author feels that it is still worthwhile to record this result and to identify the two competing notions of unstable algebras for $K$-theory.

We need to modify the above definition of an unstable $E$-cohomology algebra slightly in the case of $K$-theory by taking into account its $2$-periodicity.  Now the base point component of the $0$th space in the $\Omega$-spectrum representing $K$-theory is the classifying space $BU$ of the infinite unitary group.  Since $K^0(\text{pt}) = \bZ$ it makes sense to consider the functor
   \[
   U(-) ~=~ U^0(-) ~=~ \FRing(K(BU),-) \colon \FRing ~\to~ \Set
   \]
from the category $\FRing$ of complete Hausdorff filtered rings and continuous ring homomorphisms.  This functor can again be lifted to a comonad on the category of filtered rings. In what follows, an \emph{unstable $K$-cohomology algebra} is by definition a $U$-coalgebra for this comonad $U$.

We will define precisely what a \emph{filtered} $\lambda$-ring is below (see Definition \ref{def:filtered lambda-ring}).  This is basically a filtered ring with a $\lambda$-ring structure for which in the expression $\lambda^i(r)$, both the $\lambda$-variable and the $r$-variable are continuous.

We are now ready to state the main result of this note.

\begin{thm}
\label{thm:main}
For a complete Hausdorff filtered ring $R$, an unstable $K$-cohomology algebra structure on $R$ is exactly a filtered $\lambda$-ring structure on $R$, and vice versa.
\end{thm}

One advantage of having a result like this is that in order to study unstable algebras in the sense of Boardman-Johnson-Wilson, one has to be able to unravel the enormous amount of information encoded in a $U$-coalgebra and should compare it to more familiar structures whenever possible.  Theorem \ref{thm:main} does this for $K$-theory, identifying $U$-coalgebras with the well-studied $\lambda$-rings.

\subsection{Organization of the paper}
The rest of this paper is organized as follows.

In \S \ref{sec:filtered lambda-ring} we recall the definitions of filtered rings and $\lambda$-rings.  Then we describe the filtered ring $K(BU)$.  The main aims of the section are to define a filtered $\lambda$-ring and to observe that the completion of the $K$-theory of a space is such.

\S \ref{sec:K-cohomology} begins by recalling the notions of comonads and their coalgebras which are necessary in order to define unstable $K$-cohomology algebra.  We define the modified comonad $U$ for $K$-theory, taking into account its periodicity, and define an unstable $K$-cohomology algebra as a coalgebra over this comonad.  Particular attention is paid as to how this comonad arises from the extra structure of the space $BU$.  We also observe that the completion of the $K$-theory of a space $X$ is an unstable $K$-cohomology algebra.

In \S \ref{sec:identifying}, the final section, we then show that an unstable $K$-cohomology algebra is a filtered $\lambda$-ring and vice versa, proving Theorem \ref{thm:main}.


\section{Filtered $\lambda$-ring}
\label{sec:filtered lambda-ring}

All rings considered in this paper are assumed to be commutative, associative, and have a unit.


\subsection{Filtered ring}
\label{subsec:filtered ring}

More information about filtered objects can be found in Boardman \cite{boardman}.

A \emph{filtered ring} is a ring $R$ equipped with a directed system of ideals $I^aR$.  Directedness means that for every $I^aR$ and $I^bR$, there exists an $I^cR$ which is contained in the intersection $I^aR \cap I^bR$.  The filtration induces a filtration topology on $R$, which is Hausdorff (resp.\ complete) if and only if the natural map
   \[
   \label{eq:natural map}
   c \colon R ~\to~ \varprojlim\, R/I^aR
   \]
is injective (resp.\ surjective).  If $R$ is not already complete Hausdorff, one can always complete it by taking the inverse limit $R^\wedge = \lim\, R/I^aR$, which is complete and Hausdorff with the induced filtrations.  The ring $R^\wedge$ is called the completion of $R$.

We mention two elementary but very useful facts about the completion.  
   \begin{itemize}
   \item The map $c$ has the following universal property: If $f \colon R \to S$ is a continuous ring homomorphism to a complete Hausdorff filtered ring $S$, then $f$ factors through $c$ uniquely by a continuous ring homomorphism $f^\wedge \colon R^\wedge \to S$.  
   \item The image of $R$ in its completion is dense.
   \end{itemize}

From now on, \emph{all filtered rings will be required to be complete and Hausdorff}, unless otherwise specified.

A filtered ring homomorphism is a ring homomorphism which is continuous with respect to the filtration topology.  The category of complete Hausdorff filtered rings and continuous ring homomorphisms is denoted by $\FRing$.

When we discuss the comonad $U$ we will need to use the not-at-all obvious fact that completed tensor product is the coproduct in the category of complete Hausdorff filtered rings.  More precisely, suppose that $R = (R, \lbrace I^aR \rbrace)$ and $S = (S, \lbrace I^bS \rbrace)$ are complete Hausdorff filtered rings.  Then their \emph{completed tensor product} is defined as
   \begin{equation}
   \label{eq:completed tensor product}
   R \, \ctensor \,  S ~=~ \varprojlim_{a, \, b}\, \frac{R \otimes S}{\ker(R \otimes S \to (R/I^aR)\otimes (S/I^bS))}. 
   \end{equation}
Then $R \ctensor S$ is a complete Hausdorff filtered ring and is the coproduct of $R$ and $S$ in the category $\FRing$ (see \cite[Lemma 6.9]{boardman} for a proof).


\subsection{$\lambda$-ring}
\label{subsec:lambda-ring}

More information about $\lambda$-rings can be found in \cite{at,knutson}.

A \emph{$\lambda$-ring} is a ring $R$ equipped with functions 
   \[
   \lambda^i \colon R ~\to~ R \quad (i \geq 0)
   \]
which satisfy the following conditions.  For any integers $i, j \geq 0$ and elements $r$ and $s$ in $R$:
   \begin{itemize}
   \item $\lambda^0(r) = 1$
   \item $\lambda^1(r) = r$
   \item $\lambda^i(1) = 0$ for $i > 1$
   \item $\lambda^i(r + s) = \sum_{k = 0}^i\, \lambda^k(r)\lambda^{i-k}(s)$
   \item $\lambda^i(rs) = P_i(\lambda^1(r), \ldots , \lambda^i(r); \lambda^1(s), \ldots, \lambda^i(s))$
   \item $\lambda^i(\lambda^j(r)) = P_{i,j}(\lambda^1(r), \ldots , \lambda^{ij}(r))$.
   \end{itemize}
The polynomials $P_i$ and $P_{i,j}$ are defined as follows.  Consider the variables $\xi_1, \ldots , \xi_i$ and $\eta_1, \ldots , \eta_i$.  Denote by $s_1, \ldots , s_i$ and $\sigma_1, \ldots , \sigma_i$, respectively, the elementary symmetric functions of the $\xi$'s and the $\eta$'s.  The polynomial $P_i$ is defined by the requirement that the expression 
$P_i(s_1, \ldots , s_i; \sigma_1, \ldots , \sigma_i)$ 
be the coefficient of $t^i$ in the finite product
   \[
   \prod_{m,n=1}^i\, (1 + \xi_m \eta_n t).
   \]
Similarly, if $s_1, \ldots , s_{ij}$ are the elementary symmetric functions of $\xi_1, \ldots , \xi_{ij}$, then the polynomial $P_{i,j}$ is defined by the requirement that the expression 
$P_{i,j}(s_1, \ldots , s_{ij})$ 
be the coefficient of $t^i$ in the finite product
   \[
   \prod_{l_1 < \cdots < l_j} \, (1 + \xi_{l_1} \cdots \xi_{l_j} t).
   \]

A $\lambda$-ring map is a ring map which commutes with the $\lambda$-operations.

An example of a $\lambda$-ring is the $K$-theory $K(X)$ of a CW space $X$.  It is also a filtered ring, with filtration ideals 
   \begin{equation}
   \label{eq:profinite}
   I^a(X) ~=~ \ker(K(X) \to K(X_a)),
   \end{equation}
where $X_a$ runs through all the finite subcomplexes of $X$.  This is sometimes called the \emph{profinite filtration} on $K(X)$ (cf.\ \cite[Def.\ 4.9]{boardman}).  Although the filtration on $X$ is not unique, the isomorphism type of the resulting filtered ring structure on $K(X)$ is well-defined.

The induced map in $K$-theory of a map between spaces is, of course, a $\lambda$-ring map as well as a filtered ring map.


\subsection{The filtered ring $K(BU)$}
\label{subsec:K(BU)}

The purpose of this subsection is to describe the filtered ring $K(BU)$ in relation with the $\lambda$-operations.

Since $K$-theory is represented by the classifying space $BU$ of the infinite unitary group, the Yoneda lemma implies that $K$-theory operations can be identified with self-maps of $BU$, which in turn can be identified with the elements of $K(BU)$.  More precisely, we will consistently identify the following:
   \begin{itemize}
   \item a self-map $r \colon BU \to BU$,
   \item the corresponding element $r \in K(BU)$, and
   \item the operation $r_* \colon K(-) \to K(-)$.
   \end{itemize}

It is well-known that the $K$-theory filtered ring of $BU$ is a power series ring
   \begin{equation}
   \label{eq:K(BU)}
   K(BU) ~=~ \bZ \llbrack \lambda_1, \lambda_2, \ldots \, \rrbrack
   \end{equation}
in countably infinitely many variables $\lambda_i$ $(i > 0)$, corresponding to the operations $\lambda^i$ in $K$-theory (see, e.g., \cite[Thm.\ 4.15]{toda}).  The variable $\lambda_i$ lies in filtration exactly $2i$, and the filtered ring structure on $K(BU)$ is generated this way.


\subsection{Filtered $\lambda$-ring}
\label{subsec:filtered lambda-ring}

\begin{definition}
\label{def:filtered lambda-ring}
A \emph{filtered $\lambda$-ring} is a complete Hausdorff filtered ring $R = (R, \lbrace I^aR \rbrace)$ with a $\lambda$-ring structure such that the following two conditions hold.
   \begin{itemize}
   \item The $\lambda^i$ $(i > 0)$ is an equicontinuous family of functions.  That is, for every filtration ideal $I^aR$, there exists an $I^bR$ such that whenever $r \in I^bR$, one has $\lambda^i(r) \in I^aR$ for every $i > 0$.
   \item For every element $r \in R$ and every filtration ideal $I^aR$, there exists an integer $N > 0$ (depending on $r$ and $a$) such that whenever $\sum_{l = 1}^k \, i_l e_l \geq N$, one has $\prod_{l=1}^k\, \lambda^{i_l}(r)^{e_l} \in I^aR$.
   \end{itemize}
\end{definition}

Thus, a filtered $\lambda$-ring is essentially a filtered ring with a $\lambda$-ring structure in which the expression $\lambda^i(r)$ is continuous in both the $\lambda$- and the $r$-variables.

A filtered $\lambda$-ring map is a continuous ring homomorphism which commutes with the $\lambda$-operations.

In order that these algebraic gadgets do model the $K$-theory of spaces (completed if necessary), we have to show that $K(X)^\wedge$ is a filtered $\lambda$-ring for any CW space $X$.

\begin{prop}
\label{prop:filtered lambda-ring}
For any CW space $X$, the completion $K(X)^\wedge$ has a canonical filtered $\lambda$-ring structure for which the completion map $c \colon K(X) \to K(X)^\wedge$ is a $\lambda$-ring map.
\end{prop}

\begin{proof}
The completion $K(X)^\wedge$ is clearly a complete Hausdorff filtered ring.  Its universal property implies that the composite map
   \[
   K(X) ~\xrightarrow{\lambda^i}~ K(X) ~\xrightarrow{c}~ K(X)^\wedge
   \]
factors through $c$ uniquely via a map which we call $\lambda^i$.  Once it is shown that these $\lambda^i$ make $K(X)^\wedge$ a filtered $\lambda$-ring, it is automatically true that $c$ is a map of $\lambda$-rings.

Now since the image of $K(X)$ in its completion is dense and since $K(X)$ with its $\lambda$-operations is a $\lambda$-ring, it follows immediately that $K(X)^\wedge$ with its $\lambda^i$ is also a $\lambda$-ring.  We must still prove that these $\lambda^i$ have the required continuity properties in Definition \ref{def:filtered lambda-ring} to make $K(X)^\wedge$ a filtered $\lambda$-ring.  Again, since $K(X)$ is dense in $K(X)^\wedge$, it suffices to show that the $\lambda$-operations on $K(X)$ have these continuity properties.

To see that the $\lambda^i$ on $K(X)$ are equicontinuous, pick any filtration ideal $I^aX = \ker(K(X) \to K(X_a))$ corresponding to a finite subcomplex $X_a$ and pick any element $\alpha \in I^aX$.  Then $\alpha$ is represented by a map $\alpha \colon X \to BU$ whose restriction to $X_a$ is nullhomotopic.  If $r \in K(BU)$ is any operation at all, represented as a map $r \colon BU \to BU$, then the composite $r \circ \alpha \colon X \to BU$ is still nullhomotopic when restricted to $X_a$.  That is, the filtration ideal $I^aX$ is, in fact, closed under any $K$-theory operations, including the $\lambda$-operations.  This proves that $\lbrace \lambda^i \rbrace_{i>0}$ is an equicontinuous family of functions on $K(X)$.

To demonstrate the other continuity property, let $\alpha$ be an element of $K(X)$ and let $I^aX$ be a filtration ideal.  The element $\alpha$ is represented by a map $\alpha \colon X \to BU$, which induces a continuous ring homomorphism
   \[
   \alpha^* \colon K(BU) ~\to~ K(X)
   \]
sending an element $r \in K(BU)$ to the element $r(\alpha)$ in $K(X)$ represented by the composite
   \[
   X ~\xrightarrow{\alpha}~ BU ~\xrightarrow{r}~ BU.
   \]
That the required continuity property holds now follows from the continuity of $\alpha^*$ and the filtered ring structure on $K(BU)$ (see \S \ref{subsec:K(BU)}).
\end{proof}

\begin{cor}
\label{cor:filtered lambda-ring}
If $f \colon X \to Y$ is a map of CW spaces, then the induced map $f^{*\wedge} \colon K(Y)^\wedge \to K(X)^\wedge$ is a filtered $\lambda$-ring map.
\end{cor}


\section{Unstable $K$-cohomology algebra}
\label{sec:K-cohomology}

The main purpose of this section is to define the comonad $U$ (see eq.\ \eqref{eq:definition of U}).  Its coalgebras are by definition the unstable $K$-cohomology algebras.  We will also see that the $K$-theory of a space (completed if necessary) is a $U$-coalgebra.

We begin by recalling the concepts of comonads and their coalgebras.


\subsection{Comonad and coalgebra}
\label{subsec:comonad}
The reader can consult MacLane's book \cite{maclane} for more information on this topic.

A \emph{comonad} on a category $\bC$ is a functor $F \colon \bC \to \bC$ equipped with natural transformations $\eta \colon F \to \Id$ and $\Delta \colon F \to F^2$, called counit and comultiplication, satisfying the counital and coassociativity conditions:
   \[
   F\eta \circ \Delta ~=~ \Id_R ~=~ \eta F \circ \Delta \qquad \text{and} \qquad
   F\Delta \circ \Delta ~=~ \Delta F \circ \Delta.
   \]
The natural transformations $\Delta$ and $\eta$ are often omitted from the notation and we speak of $F$ as a comonad.

If $F$ is a comonad on a category $\bC$, then an $F$-{\it coalgebra structure} on an object $X$ of $\bC$ is a morphism $\xi \colon X \to FX$ in $\bC$, called the structure map, satisfying the counital and coassociativity conditions:
   \begin{equation}
   \label{eq:counit}
   F\xi \circ \xi ~=~ \Delta \circ \xi \qquad \text{and} \qquad \eta \circ \xi ~=~ \Id_X.
   \end{equation}
We sometimes abuse notation and say that $X$ is an $F$-coalgebra, leaving the structure map $\xi$ implicit.

A map of $F$-coalgebras $g \colon (X, \xi_X) \to (Y, \xi_Y)$ consists of a morphism $g \colon X \to Y$ in $\bC$ such that $\xi_Y \circ g = Fg \circ \xi_X$.


\subsection{The comonad $U$ for $K$-theory}
\label{subsec:comonad U}

The discussion in this section follows closely \S 8 of Boardman, Johnson, and Wilson \cite{bjw}, except that we take into account the periodicity of $K$-theory and consider only the degree $0$ part.  We will first define the comonad $U$ and then discuss its ring structure (when applied to a filtered ring), filtration, comultiplication and counit.  After that we will define unstable $K$-cohomology algebra and observe that the argument in Boardman-Johnson-Wilson shows that the $K$-theory of any space (completed if necessary) is such.

We will use the Yoneda lemma many times without explicitly mentioning it.

\subsubsection{Definition of $U$}
\label{subsubsec:definition of U}

Let, then, $R = (R, \lbrace I^aR \rbrace)$ be an arbitrary complete Hausdorff filtered ring.  Define the functor $U$ (to $\Set$ only at the moment) to be
   \begin{equation}
   \label{eq:definition of U}
   UR ~=~ \FRing(K(BU), R),
   \end{equation}
the set of continuous ring homomorphisms from $K(BU)$ to $R$.

\subsubsection{Ring structure on $UR$}
\label{subsubsec:UR}

Let $X$ be an arbitrary CW space.

There are maps 
   \[
   \label{eq:mu and phi}
   \mu,\, \phi \colon BU \times BU ~\rightrightarrows~ BU
   \]
which induce the natural addition (by $\mu$) and multiplication (by $\phi$) structure on the $K$-theory $K(X)$ of a space $X$.  These maps satisfy certain associativity, commutativity, etc.\ conditions which make $K(X)$ a commutative ring with unit.  Thanks to the Kunneth homeomorphism $K(BU \times BU) \cong K(BU) \ctensor K(BU)$ \cite[Thm.\ 4.19]{boardman}, they induce in $K$-theory the following maps:
   \[
   \label{eq:mu and phi in K}
   \mu^*,\, \phi^* \colon \bZ \llbrack \lambda_1, \lambda_2, \ldots \, \rrbrack 
   ~\rightrightarrows~ \bZ \llbrack \lambda_1, \lambda_2, \ldots \, \rrbrack ~\ctensor~ 
   \bZ \llbrack \lambda_1, \lambda_2, \ldots \, \rrbrack.
   \]
Therefore, for each integer $k \geq 1$, we can write
   \[
   \mu^*(\lambda_k) ~=~ \sum_{\alpha}\, r_\alpha^\prime \, \otimes r_\alpha^{\prime\prime}
   \]
for some elements $r_\alpha^\prime$ and $r_\alpha^{\prime\prime}$ in $K(BU)$.  By precomposition this becomes the Cartan formula for a sum
   \begin{equation}
   \label{eq:Cartan for sum}
   \begin{split}
   \lambda^k(x + y) 
   &~=~ \sum_{\alpha}\, r_\alpha^\prime(x) r_\alpha^{\prime\prime}(y) \quad (x, y \in K(X)) \\
   &~=~ \sum_{i + j = k}\, \lambda^i(x) \lambda^j(y)
   \end{split}
   \end{equation}
where we have used the usual convention $\lambda^0(x) = 1$.   Since eq.\ \eqref{eq:Cartan for sum} holds for an arbitrary space $X$ and any elements $x$ and $y$ in $K(X)$, we must have that
   \begin{equation}
   \label{eq:mu}
   \mu^*(\lambda_k) ~=~ \sum_{i + j = k}\, \lambda_i \otimes \lambda_j.
   \end{equation}

A similar reasoning, using the property 
   \[
   \lambda^k(xy) ~=~ P_k(\lambda^1x, \ldots , \lambda^kx; \lambda^1y, \ldots , \lambda^ky),
   \]
leads to the formula
   \begin{equation}
   \label{eq:phi}
   \phi^*(\lambda_k) ~=~ P_k(\lambda_1 \otimes 1, \ldots , \lambda_k \otimes 1; 1 \otimes \lambda_1, \ldots , 1 \otimes \lambda_k).
   \end{equation}

Now suppose that $f$ and $g$ are elements of $UR = \FRing(K(BU),R)$.   Their sum and product, $f+g$ and $fg$, both have the form   
   \[
   K(BU) ~\to~ K(BU) \ctensor K(BU) ~\xrightarrow{f \ctensor g}~ R \ctensor R ~\xrightarrow{\text{multiplication}}~ R,
   \]
where for $f+g$ (resp.\ $fg$) the left-hand map is $\mu^*$ (resp.\ $\phi^*$).  
Then eq. \eqref{eq:mu} and \eqref{eq:phi} imply that on the element $\lambda_k \in K(BU)$ these maps can be expressed as the Cartan formulas
   \begin{equation}
   \label{eq:sum and product}
   \begin{split}
   (f + g)(\lambda_k) &~=~ \sum_{i+j = k}\, f(\lambda_i) g(\lambda_j) \\
   (fg)(\lambda_k)    &~=~ P_k(f(\lambda_1), \ldots , f(\lambda_k); g(\lambda_1), \ldots , g(\lambda_k)).
   \end{split}
   \end{equation}
It follows that the additive and multiplicative identities, $0_{UR}$ and $1_{UR}$, of $UR$ are given by the maps
   \begin{equation}
   \label{eq:0 and 1}
   \begin{split}
   0_{UR} & \colon \lambda_k \mapsto 0 \quad (k > 0) \\
   1_{UR} & \colon \lambda_k \mapsto \begin{cases} 1 \text{ if } k = 1 \\
                                    0 \text{ if } k > 1.
                      \end{cases}
   \end{split}
   \end{equation}
The second equality follows from the fact that $P_k(x_1, \ldots , x_k; 1, 0, \ldots , 0) = x_k$.

\subsubsection{Filtration on $UR$}
\label{subsubsec:filtration on UR}

The ring $UR = \FRing (K(BU),R)$ is filtered by the ideals
   \begin{equation}
   \label{eq:I^aUR}
   \begin{split}
   I^aUR &~=~ \ker(UR \to U(R/I^aR)) \\
         &~=~ \lbrace f \in UR \, \colon \, f(\lambda_k) \in I^aR \text{ for all } k>0 \rbrace.
   \end{split}
   \end{equation}
It is easy to see that the surjective map $UR \to U(R/I^aR)$ has kernel $I^aUR$.  Thus, since $R$ is complete Hausdorff, it follows that
   \[
   \begin{split}
   UR &~=~ \FRing(K(BU),R) \\
      &~=~ \lim_a \, \FRing(K(BU), R/I^aR) \\
      &~=~ \lim_a \, (UR)/I^aUR
   \end{split}
   \]
That is, $UR$ is also a complete Hausdorff filtered ring.  Note that the indexing set for the filtration of $UR$ is the same as that for $R$.

So far we have seen that $U$ is a functor on the category $\FRing$ of complete Hausdorff filtered rings.

\subsubsection{Comultiplication and counit}
\label{subsubsec:comultiplication}

In order to make $U$ a comonad on $\FRing$, we still need the natural transformations $\Delta \colon U \to U^2$ and $\eta \colon U \to \Id$.  Let us begin with the former.  

There is a filtered ring map
   \[
   \begin{split}
   \rho \colon K(BU) &~\xrightarrow~ U(K(BU)) \\
                   g &~\mapsto~ (f \mapsto f \circ g)
   \end{split}
   \]
Using the formula $\lambda^i\lambda^j(x) = P_{i,j}(\lambda^1x, \ldots , \lambda^{ij}x)$, we can express the map $\rho$ in terms of the elements $\lambda_i \in K(BU)$ as follows:
   \begin{equation}
   \label{eq:composition}
   \begin{split}
   \rho(\lambda_j)(\lambda_i) &~=~ \lambda^i \circ \lambda^j \\
           &~=~ P_{i,j}(\lambda_1, \ldots , \lambda_{ij}).
   \end{split}
   \end{equation}
Here $\lambda^i \circ \lambda^j$ is the element in $K(BU)$ represented by the composition of the $K$-theory operations $\lambda^i$ and $\lambda^j$.  
Now if $f$ is an element in $UR$, then its image under the comultiplication map $\Delta_R \colon UR ~\to~ U^2R$ is the composite map
   \begin{equation}
   \label{eq:comultiplication for U}
   \Delta_R f ~=~ (Uf) \circ \rho \colon K(BU) ~\to~ U(K(BU)) ~\to~ UR.
   \end{equation}

As for the counit $\eta \colon U \to \Id$, it is defined by
   \begin{equation}
   \label{eq:counit for U}
   \eta_Rf ~=~ f(\lambda_1);
   \end{equation}
that is, $\eta_R$ is simply the evaluation map at $\lambda_1$.

\subsubsection{Unstable $K$-cohomology algebra}

The proof in Boardman-Johnson-Wilson \cite[Thm.\ 8.8(a)]{bjw} that their $U$ is a comonad on the category of complete Hausdorff filtered $E^*$-algebras carries over almost without change to show the following.

\begin{prop}
\label{prop:U is a comonad}
The functor $U$ defined in eq.\ \eqref{eq:definition of U} together with $\Delta$ and $\eta$ above is a comonad on the category $\FRing$ of complete Hausdorff filtered rings.
\end{prop}

Following Boardman, Johnson, and Wilson we now make the following definition.

\begin{definition}
\label{def:unstable K-cohomology algebra}
An \emph{unstable $K$-cohomology algebra} is a $U$-coalgebra for the comonad $U$ in Proposition \ref{prop:U is a comonad}.  A map of unstable $K$-cohomology algebras is a map of $U$-coalgebras.
\end{definition}

Now given any CW space $X$, composition of maps yields a continuous map
   \[
   \circ \, \colon K(BU) \times K(X) ~\to~ K(X)
   \]
which gives, after completion, the map
   \[
   K(BU) \times K(X)^\wedge ~\to~ K(X)^\wedge.
   \]
Taking its adjoint we obtain a map
   \begin{equation}
   \label{eq:rho_X}
   \rho_X \colon K(X)^\wedge ~\to~ \FRing(K(BU),K(X)^\wedge) ~=~ U(K(X)^\wedge).
   \end{equation}
The argument of \cite[Thm.\ 8.11(a)]{bjw} in Boardman-Johnson-Wilson, which shows that their analogous map $\rho_X \colon E^*(X)^\wedge \to U(E^*(X)^\wedge)$ makes $E^*(X)^\wedge$ a $U$-coalgebra, now gives

\begin{prop}
\label{prop:K(X) is a U-coalgebra}
The map $\rho_X$ in eq.\ \eqref{eq:rho_X} makes $K(X)^\wedge$ an unstable $K$-cohomology algebra. 
\end{prop}

This proposition is also a consequence of Proposition \ref{prop:filtered lambda-ring} and Theorem \ref{thm:main} in the Introduction.

\begin{cor}
\label{cor:U-coalgebra}
If $f \colon X \to Y$ is a map of CW spaces, then the induced map $f^{*\wedge} \colon K(Y)^\wedge \to K(X)^\wedge$ is a map of unstable $K$-cohomology algebras.
\end{cor}


\section{Identifying unstable $K$-cohomology algebra with filtered $\lambda$-ring}
\label{sec:identifying}

In this final section we will prove Theorem \ref{thm:main} in the Introduction.

So let $R = (R, \lbrace I^aR \rbrace)$ be an arbitrary complete Hausdorff filtered ring.  Suppose that $R$ has the structure of an unstable $K$-cohomology algebra, $\xi \colon R \to UR$.  We must show that this gives a filtered $\lambda$-ring structure on $R$.  Recall that
    \[
    UR ~=~ \FRing(K(BU),R) ~=~ \FRing(\bZ \llbrack \lambda_1, \lambda_2, \ldots \, \rrbrack, R).
    \]

We define the operations 
   \[
   \label{eq:lambda on R}
   \lambda^i \colon R ~\to~ R \quad (i \geq 0)
   \]
by setting $\lambda^0 \equiv 1$ and for $i > 0$, 
   \begin{equation}
   \label{eq:lambda^i}
   \lambda^i(r) ~\buildrel \text{def} \over =~ (\xi r)(\lambda_i) \quad (r \in R).
   \end{equation}
We claim that these operations make $R$ into a filtered $\lambda$-ring; that is, a $\lambda$-ring structure on $R$ together with the two continuity properties in Definition \ref{def:filtered lambda-ring}.  The argument is divided into six steps, the first one for the continuity properties and the rest for the $\lambda$-ring structure.

\bigskip
\emph{Step 1}. We first check the continuity properties in Definition \ref{def:filtered lambda-ring}.  To see that the family $\lbrace \lambda^i \rbrace_{i>0}$ is equicontinuous, let $I^aR$ be a filtration ideal.  We must show that there exists an $I^bR$ such that $\lambda^i(I^bR) \subset I^aR$ for every $i > 0$.  Since $\xi$ is continuous, given $I^aUR$, there exists $I^bR$ such that $\xi(I^bR) \subset I^aUR$.  That is, if $r \in I^bR$, then $\lambda^i(r) = (\xi r)(\lambda_i) \in I^aR$ for every $i > 0$.  This shows that $\lbrace \lambda^i \rbrace_{i>0}$ is an equicontinuous family of functions on $R$.

To check the second continuity property in Definition \ref{def:filtered lambda-ring}, let $r$ be an element in $R$ and let $I^aR$ be a filtration ideal.  The element $\xi r \in UR$ is a continuous ring homomorphism from $K(BU)$ to $R$.  Thus, given $I^aR$, there exists an integer $N > 0$ such that whenever $\alpha \in K(BU)$ has filtration strictly great than $N$, then $(\xi r)(\alpha) \in I^aR$.  Now if $\sum_{l=1}^k\, i_l e_l \geq N$, then the filtration of the element $\prod_{l=1}^k\, \lambda_{i_l}^{e_l} \in K(BU)$ is $\sum_{l=1}^k\, 2i_l e_l > N$, and so we have
   \[
   \prod_{l=1}^k\, \lambda^{i_l}(r)^{e_l} ~=~ (\xi r)\left( \prod_{l=1}^k \, \lambda_{i_l}^{e_l} \right) \in I^aR.
   \]
This proves the desired continuity property.

\bigskip
\emph{Step 2}.  We check that $\lambda^1$ is the identity map on $R$.  Recall that the counit $\eta_R \colon UR \to R$ is the evaluation map at $\lambda_1$.  Since there is an equality $\Id_R = \eta_R \xi$, it follows that for any element $r$ in $R$, we have that
   \[
   \begin{split}
   \lambda^1(r) &~=~ (\xi r)(\lambda_1) \\
                &~=~ (\eta_R \xi)(r) \\
                &~=~ r.
   \end{split}
   \]
So $\lambda^1$ is the identity map on $R$.

\bigskip
\emph{Step 3}.  We check that $\lambda^i(1) = 0$ for any $i > 1$.  Denoting the multiplicative identity of $UR$ by $1_{UR}$ (see eq.\ \eqref{eq:0 and 1}), we have that for $i > 1$,
   \[
   \begin{split}
   \lambda^i(1) &~=~ (\xi 1)(\lambda_i) \\
                &~=~ 1_{UR}(\lambda_i)  \\
                &~=~ 0
   \end{split}
   \]
as desired.

\bigskip
\emph{Step 4}.  Now we show the Cartan formula for a sum of any two elements $x$ and $y$ in $R$.  Using the additivity of $\xi$ and eq.\ \eqref{eq:sum and product}, we calculate
   \[
   \begin{split}
   \lambda^k(x + y) &~=~ \xi(x + y)(\lambda_k) \\
                    &~=~ (\xi x + \xi y)(\lambda_k) \\
                    &~=~ \sum_{i + j=k}\, \lbrace (\xi x)(\lambda_i)\rbrace 
                                          \lbrace (\xi y)(\lambda_j)\rbrace \\
                    &~=~ \sum_{i + j=k}\, \lambda^i(x)\lambda^j(y)
   \end{split}
   \]
This proves the Cartan formula for a sum.

\bigskip
\emph{Step 5}.  The Cartan formula for $\lambda^k(xy)$ is proved similarly, using the multiplicativity of $\xi$ and eq.\ \eqref{eq:sum and product}.

\bigskip
\emph{Step 6}.  Finally, we show that $\lambda^i\lambda^j(x) = P_{i,j}(\lambda^1x, \ldots, \lambda^{ij}x)$.   Let $x$ be an element in $R$.  Then, using the coassociativity of $\xi$ (see eq.\ \eqref{eq:counit}) and eq.\ \eqref{eq:composition} and \eqref{eq:comultiplication for U}, we calculate
   \[
   \begin{split}
   \lambda^i \lambda^j (x) 
   &~=~ \lambda^i(\xi x(\lambda_j)) \\
   &~=~ \lbrace \xi (\xi x(\lambda_j))\rbrace (\lambda_i) \\
   &~=~ \lbrace\Delta_R(\xi x)(\lambda_j) \rbrace (\lambda_i) \\
   &~=~ \lbrace U(\xi x) \circ \rho \rbrace (\lambda_j)(\lambda_i) \\
   &~=~ (\xi x) P_{i,j}(\lambda_1, \ldots , \lambda_{ij}) \\
   &~=~ P_{i,j}(\lambda^1 x, \ldots , \lambda^{ij}x)
   \end{split}
   \]
as desired.

We have shown that the operations $\lambda^i$ in eq.\ \eqref{eq:lambda^i} make $R$ into a $\lambda$-ring which also satisfies the two continuity properties in Definition \ref{def:filtered lambda-ring}.  Therefore, the unstable $K$-cohomology algebra structure on $R$ gives a filtered $\lambda$-ring structure on $R$.  This proves half of Theorem \ref{thm:main}.

The above argument can easily be reversed to show that a filtered $\lambda$-ring structure on $R$ yields, via eq.\ \eqref{eq:lambda^i}, an unstable $K$-cohomology algebra structure on $R$.  Indeed, the two continuity properties in the definition of a filtered $\lambda$-ring make sure that both the proposed structure map $\xi \colon R \to UR$ and $\xi(r) \colon K(BU) \to R$ for any $r \in R$ are continuous.  The Cartan formulas for $\lambda^n(x+y)$ and $\lambda^n(xy)$ imply the additivity and multiplicativity, respectively, of the structure map $\xi$, and the property about $\lambda^i\lambda^j(x)$ leads to the coassociativity of $\xi$.  That $\xi$ is counital follows from the condition that $\lambda^1$ is the identity on $R$.

The proof of Theorem \ref{thm:main} is complete.


\end{document}